\documentclass[final,leqno]{article}
\usepackage{epsfig,psfrag,amsmath,amssymb}
\newtheorem{mytheorem}{Theorem}
\newtheorem{mylemma}{Lemma}

\newtheorem{myremark}[mytheorem]{Remark}

\newtheorem{mycorollary}{Corollary}

\newcommand{\N}{{\sf N\hspace*{-1.0ex}\rule{0.15ex}{1.3ex}\hspace*{1.0ex}}}

\newcommand{\BA}{\begin{eqnarray}}
\newcommand{\EA}{\end{eqnarray}}
\newcommand{\BE}{\begin{equation}}
\newcommand{\EE}{\end{equation}}
\newcommand{\ba}{\begin{array}}
\newcommand{\ea}{\end{array}}
\newcommand{\baa}{\begin{eqnarray*}}
\newcommand{\eaa}{\end{eqnarray*}}
\newcommand{\be}{\begin{equation}}
\newcommand{\ee}{\end{equation}}

\def\N{{\cal N}}

\def\u{{\bf u}}
\def\w{{\bf w}}

\def\Rtilde{{\widetilde{R}}}
\def\Btilde{{\widetilde{B}}}

\def\Atilde{{\widetilde{A}}}
\def\Wtilde{{\widetilde{W}}}

\def\>{\raisebox{-1ex}{$\; \stackrel{\textstyle >}{\sim } \; $}}
\def\<{\raisebox{-1ex}{$ \; \stackrel{\textstyle <} {\sim } \; $}}
\newcommand{\vvec}[1]{{\mathbf{#1}}}

\newcommand{\EQ}[1]{(\ref{equation:#1})}
\newcommand{\LA}[1]{\ref{lemma:#1}}

\def\beginproof{\indent {\it Proof:~}}

\title{A non-overlapping domain decomposition method for incompressible Stokes equations with continuous pressure}

\author{Jing Li\thanks{Department of Mathematical Sciences, Kent State University, Kent, OH 44242, {\tt li@math.kent.edu}, {\tt http://www.math.kent.edu/$\sim$li/}.} \and Xuemin Tu\thanks{Department of Mathematics, University of Kansas, 1460 Jayhawk Blvd, Lawrence, KS 66045-7594,  {\tt xtu@math.ku.edu}, {\tt http://www.math.ku.edu/$\sim$xtu/}. This author's work was supported in part by National Science Foundation contract DMS-1115759.} }

\begin{document}

\maketitle
\pagestyle{myheadings} \thispagestyle{plain} \markboth{JING LI AND XUEMIN TU}{DOMAIN DECOMPOSITION FOR INCOMPRESSIBLE STOKES}

\begin{abstract}
A non-overlapping domain decomposition algorithm is proposed to solve
the linear system arising from mixed finite element approximation of incompressible Stokes equations. A continuous finite element space for the pressure is used. In the proposed algorithm, Lagrange multipliers are used to enforce continuity of the velocity component across the subdomain domain boundary. The continuity of the pressure component is enforced in the primal form, i.e., neighboring subdomains share the same pressure degrees of freedom on the subdomain interface and no Lagrange multipliers are needed. After eliminating all velocity variables and the independent subdomain interior parts of the pressures, a symmetric positive semi-definite linear system for the subdomain boundary pressures and the Lagrange multipliers is formed and solved by a preconditioned conjugate gradient method. A lumped preconditioner is studied and the condition number bound of the preconditioned operator is proved to be independent of the number of subdomains for fixed subdomain problem size. Numerical experiments demonstrate the convergence rate of the proposed algorithm.
\end{abstract}

{\bf keywords}
domain decomposition, incompressible Stokes, FETI-DP, BDDC

{\bf AMS}
65F10, 65N30, 65N55

\section{Introduction}

Domain decomposition methods have been studied well for solving
incompressible Stokes equations and similar saddle-point problems; see, e.g., \cite{Kla98, 8pavwid, li05, paulo2003, Doh04, li06, kim06, Tu:2005:BPP, Tu:2005:BPD, LucaOlofStef}. In many of those work, special
care need be taken to deal with the divergence-free constraints across
subdomain boundaries, which often lead to large coarse level
problems. The large coarse level problem will be a bottleneck in
large scale parallel computations, and additional efforts in the algorithm are needed to reduce its impact, cf.~\cite{Tu:2004:TLB,Tu:2005:TLB,Tudd16,Klawonn:2005:IFM,Dohrmann:2005:ABP,KimTu,Tu:2011:TLBS}. Some recent progress has been made by Dohrmann and Widlund~\cite{Doh09, Doh10} for the almost incompressible elasticity, where the
coarse level space is built from discrete subdomain saddle-point harmonic extensions of certain subdomain interface cut-off functions and its dimension is much smaller than those in the previous studies. Kim and Lee~\cite[with Park]{kim11, kim102, kim10} studied both the FETI-DP and BDDC algorithms for incompressible Stokes equations where a lumped preconditioner is used and reduction in the dimension of the coarse level space is also achieved.

In most above mentioned applications and analysis of domain
decomposition methods for incompressible Stokes equations, the mixed
finite element space contains discontinuous pressures. Application of
discontinuous pressures in domain decomposition methods is
natural. The decomposing of the pressure components to independent
subdomains can be handled conveniently and no continuity of pressures
across the subdomain boundary need be enforced. However, a big class
of mixed finite elements used for solving incompressible Stokes and
Navier-Stokes equations have continuous pressures, e.g., the well
known Taylor-Hood type~\cite{Taylor}. There have been a variety of
approaches using continuous pressures in domain decomposition methods
for solving incompressible Stokes equations, e.g., by Goldfeld \cite{pauloPhD}, by \v{S}\'istek {\em et. al.} \cite{sis11}, and by Benhassine and Bendali \cite{ben10}. In their work, an indefinite system of linear equations need be solved, either by a generalized minimal residual method or simply by a conjugate gradient method. To the best of our knowledge, no scalable convergence rate has been proved analytically for any of those approaches using continuous pressures.

In this paper, we propose a non-overlapping domain decomposition
algorithm for solving incompressible Stokes equations with continuous
pressure finite element space. The scalability of its
convergence rate is proved.
In this algorithm, the subdomain boundary velocities are dealt with in
the same way as in the FETI-DP method: a few for each subdomain are
selected as the coarse level primal variables, which
are shared by neighboring subdomains; the others are subdomain
independent and Lagrange multipliers are used to enforce their
continuity. The subdomain boundary pressure degrees of freedom are all
in the primal form. They are shared by neighboring subdomains and no
Lagrange multipliers are needed for their continuity. After
eliminating all velocity variables and the independent subdomain
interior parts of the pressures, the system for the subdomain boundary
pressures and the Lagrange multipliers is shown to be symmetric
positive semi-definite. A preconditioned conjugate gradient method
with  a lumped preconditioner is studied. As strong condition number bounds as for the scalar elliptic case are established. In the proposed algorithm and in the estimate of its condition number bound, no additional coarse level variables, except those necessary for solving scalar elliptic problems, are required for incompressible Stokes problems. The resulting coarse level problem is also symmetric positive definite.

To stay focused on the purpose of this paper, the discussion of the
proposed algorithm and its analysis are based on two-dimensional
problems, even though the same approach can be extended to the three-dimensional case without substantial obstacles. It is also worth pointing out that the domain decomposition algorithm and its analysis presented in this paper apply equally well, with only minor modifications, to the case where discontinuous pressures are used in the mixed finite element space.

The remainder of this paper is organized as follows. The finite element discretization of the incompressible Stokes equation is introduced in Section \ref{section:FEM}. A domain decomposition approach is described in Section~\ref{section:DDM}. The system for the subdomain boundary pressures and the Lagrange multipliers is derived in Section~\ref{section:Gmatrix}. Section \ref{section:techniques} provides some techniques used in the condition number bound estimate. In Section~\ref{section:lumped}, a lumped preconditioner is proposed and a scalable condition number bound of the preconditioned operator is established. At the end, in Section~\ref{section:numerics}, numerical results for solving a two-dimensional incompressible Stokes problem are shown to demonstrate the convergence rate of the proposed algorithm.

\section{Finite element discretization}
\label{section:FEM}

We consider solving the following incompressible Stokes problem on a
bounded, two-dimensional polygonal domain $\Omega$  with a
Dirichlet boundary condition,
\begin{equation}
\label{equation:Stokes}
\left\{
\begin{array}{rcll}
-\Delta {\bf u} + \nabla p & = & {\bf f}, & \mbox{ in } \Omega \mbox{ , } \\
-\nabla \cdot {\bf u}       & = & 0, & \mbox{ in } \Omega \mbox{ , } \\
{\bf u}                    & = & {\bf u}_{\partial \Omega}, & \mbox{ on } \partial \Omega \mbox{ , }\\
\end{array}\right.
\end{equation}
where the boundary data ${\bf u}_{\partial \Omega}$ satisfies the compatibility
condition $\int_{\partial \Omega} {\bf u}_{\partial \Omega} \cdot {\bf n} = 0$. For
simplicity, we assume that ${\bf u}_{\partial \Omega} = {\bf 0}$ without losing any
generality.

The weak solution of \EQ{Stokes} is given by: find $\vvec{u} \in
\left(H^1_0(\Omega)\right)^2 = \{ \vvec{v} \in (H^1(\Omega))^2 ~
\big| ~ \vvec{v} = \vvec{0} \mbox{ on }
\partial \Omega \}$ and $p \in
L^2(\Omega)$, such that
\begin{equation}
\label{equation:bilinear} \left\{
\begin{array}{lcll}
a(\vvec{u}, \vvec{v}) + b(\vvec{v}, p) & = & (\vvec{f}, \vvec{v}),
& \forall \vvec{v}\in \left(H^1_0(\Omega)\right)^2 , \\ [0.5ex]
b(\vvec{u}, q) & = & 0, & \forall q \in L^2(\Omega) \mbox{ , }
\\
\end{array} \right.
\end{equation}
where
\[
a(\vvec{u}, \vvec{v})= \int_{\Omega} \nabla{\bf u} \cdot \nabla{\bf v}, \quad
b(\vvec{u},q) = -\int_{\Omega} (\nabla \cdot \vvec{u}) q, \quad
(\vvec{f}, \vvec{v}) = \int_{\Omega} \vvec{f} \cdot \vvec{v}.
\]
We note that the solution of \EQ{bilinear} is not unique, with the pressure $p$ different up to an additive constant.

\begin{figure}[h]
\begin{center}
\includegraphics[scale=.95]{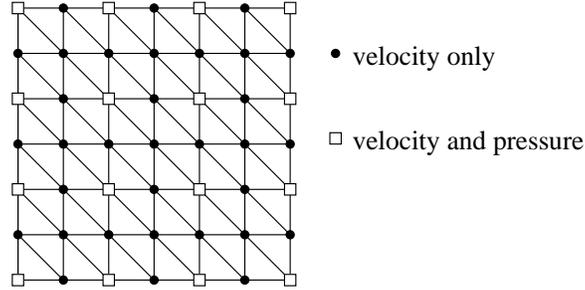}
\end{center}
\caption{\label{figure:TaylorHood} A modified Taylor-Hood mixed finite element}
\end{figure}

A modified Taylor-Hood mixed finite element is used in this paper to solve \EQ{bilinear}. The domain $\Omega$ is triangulated into shape-regular elements of characteristic size $h$. The pressure finite element space, $Q \subset L^2(\Omega)$, is taken as the space of continuous piecewise linear functions on the triangulation. The velocity finite element space,  $\vvec{W} \in \left(H^1_0(\Omega)\right)^2$, is formed by the continuous piecewise linear functions on the finer triangulation obtained by dividing each triangle into four subtriangles by connecting the middle points of its edges. A demonstration of this mixed finite element on a triangulation of a square domain is shown in Figure \ref{figure:TaylorHood}.

The finite element solution $(\vvec{u}, p) \in \vvec{W} \bigoplus Q$ of \EQ{bilinear} satisfies
\begin{equation}
\label{equation:matrix} \left[
\begin{array}{cccc}
A     &  B^T \\
B     &  0   \\
\end{array}
\right] \left[
\begin{array}{c}
{\bf u}     \\
p   \\
\end{array}
\right] = \left[
\begin{array}{l}
{\bf f}        \\
0      \\
\end{array}
\right] ,
\end{equation}
where $A$, $B$, and $\vvec{f}$  represent respectively the restrictions of $a( \cdot , \cdot )$, $b(\cdot, \cdot )$ and
$(\vvec{f} , \cdot)$ to the finite-dimensional spaces $\vvec{W}$ and $Q$. We use the same notation in this paper to represent both a finite element function and the vector of its nodal values.

The coefficient matrix in \EQ{matrix} is rank deficient. $A$ is symmetric positive definite. The kernel of $B^T$, denoted by $Ker(B^T)$, is the space of all constant pressures in $Q$. The range of $B$, denoted by $Im(B)$,  is orthogonal to $Ker(B^T)$ and is the subspace of $Q$ consisting of all vectors with zero average. The solution of \EQ{matrix} always exists and is uniquely determined when the pressure is considered in the quotient space $Q/Ker(B^T)$. In this paper, when $q \in Q/Ker(B^T)$, $q$ always has zero average. For a more general right-hand side vector $({\bf f}, ~ g)$ given in \EQ{matrix}, the existence of its solution requires that $g \in Im(B)$, i.e., $g$ has zero average.

The modified Taylor-Hood mixed finite element
space $\vvec{W} \times Q$, as shown in Figure \ref{figure:TaylorHood}, is inf-sup stable in the sense that there exists
a positive constant $\beta$, independent of $h$, such that
\begin{equation}
\label{equation:infsup} \sup_{\vvec{w} \in \vvec{W}}
\frac{b(\vvec{w},q)}{|\vvec{w}|_{H^1}} \geq \beta \|q\|_{L^2},
\hspace{0.5cm} \forall q \in Q/Ker(B^T),
\end{equation}
cf.~\cite[Chapter III, \S 7]{braess},  or equivalently in matrix/vector form,
\begin{equation}
\label{equation:infsupMatrix} \sup_{{\bf w} \in {\bf W}}
\frac{\left< q, B \vvec{w} \right>^2}{\left< \vvec{w},  A \vvec{w} \right>} \geq \beta^2 \left< q, Z q \right>,
\hspace{0.5cm} \forall q \in Q/Ker(B^T).
\end{equation}
Here, as always in this paper, $\left< \cdot, \cdot \right>$ represents the inner product of two vectors.
The matrix $Z$ represents the mass matrix defined on the pressure
finite element space $Q$, i.e., for any $q \in Q$, $\|q\|_{L^2}^2 =
\left< q, Z q \right>$. It is easy to see, cf.~\cite[Lemma B.31]{Toselli:2004:DDM}, that $Z$ is spectrally equivalent to $h^2 I$ for two-dimensional problems, where $I$ represents the identity matrix of the same dimension, i.e., there exist positive constants $c$ and $C$, such that
\be
\label{equation:massmatrix}
c h^2 I \leq Z \leq C h^2 I.
\ee
Here, as in other places of this paper, $c$ and $C$  represent generic positive constants which are independent of the mesh size $h$ and the subdomain diameter $H$ (discussed in the following section).


\section{A non-overlapping domain decomposition approach}
\label{section:DDM}

The domain $\Omega$ is decomposed into $N$ non-overlapping polygonal subdomains $\Omega_i$, $i = 1, 2, ..., N$. Each subdomain is the union of a bounded number of elements, with the diameter of the subdomain in the order of $H$. The nodes on the interface of
neighboring subdomains match across the subdomain boundaries $\Gamma = {(\cup\partial\Omega_i)} \backslash
\partial\Omega$. $\Gamma$ is composed of subdomain edges, which are regarded as open subsets of
$\Gamma$, and of the subdomain vertices, which are end points of edges.

The velocity and pressure finite element spaces ${\bf W}$ and $Q$ are decomposed into
\[
{\bf W} = {\bf W}_I \bigoplus {\bf W}_{\Gamma}, \quad
Q = Q_I \bigoplus Q_\Gamma,
\]
where ${\bf W}_I$ and $Q_I$ are direct sums of independent subdomain interior velocity spaces ${\bf W}^{(i)}_I$, and interior pressure spaces $Q^{(i)}_I$, respectively, i.e.,
$$
{\bf W}_I = \bigoplus_{i=1}^{N}{\bf W}^{(i)}_I, \quad Q_I =
\bigoplus_{i=1}^{N}Q^{(i)}_I.
$$
${\bf W}_{\Gamma}$ and $Q_\Gamma$ are subdomain boundary velocity and pressure spaces, respectively. All
functions in ${\bf W}_{\Gamma}$  and $Q_\Gamma$ are continuous across the subdomain boundaries $\Gamma$; their degrees of freedom are shared by neighboring subdomains.

To formulate our domain decomposition algorithm, we introduce a partially sub-assembled subdomain boundary
velocity space $\vvec{\Wtilde}_{\Gamma}$,
\[
\vvec{\Wtilde}_{\Gamma} = \vvec{W}_{\Pi} \bigoplus
\vvec{W}_{\Delta} = \vvec{W}_{\Pi} \bigoplus \left(
\bigoplus_{i=1}^N \vvec{W}^{(i)}_\Delta \right).
\]
Here, $\vvec{W}_{\Pi}$ is the continuous, coarse level, primal velocity space which is typically spanned by subdomain
vertex nodal basis functions, and/or by interface edge basis functions with constant values, or with values of positive
weights on these edges. The primal, coarse level velocity degrees of freedom are shared by neighboring subdomains. The complimentary space $\vvec{W}_{\Delta}$ is the direct sum of independent subdomain dual interface velocity spaces $\vvec{W}_{\Delta}^{(i)}$, which correspond to the remaining subdomain boundary velocity degrees of freedom and are spanned by basis functions which vanish at the primal degrees of freedom. Thus, an element in the space $\vvec{\Wtilde}_{\Gamma}$ typically has a continuous primal velocity component and a discontinuous dual velocity component.

The functions ${\bf w}_{\Delta}$ in ${\bf W}_{\Delta}$ are in general not continuous across $\Gamma$. To enforce their continuity, we define a boolean matrix $B_\Delta$  constructed from $\{0,1,-1\}$. On each row of $B_\Delta$, there are only two non-zero entries, $1$ and $-1$, corresponding to the same velocity degree of freedom on each subdomain boundary node, but attributed to two neighboring subdomains, such that for any ${\bf w}_{\Delta}$ in ${\bf W}_{\Delta}$, each row of $B_\Delta {\bf w}_{\Delta} = 0$ implies that these two degrees of freedom from the two neighboring subdomains be the same. When non-redundant continuity constraints are enforced, $B_\Delta$ has full row rank. We denote the range of $B_\Delta$ applied on ${\bf W}_{\Delta}$ by $\Lambda$, the vector space of the Lagrange multipliers.

In order to define a certain subdomain boundary scaling operator, we introduce a positive scaling factor $\delta^{\dagger}(x)$ for each node $x$ on the subdomain boundary $\Gamma$. Let $\N_x$ be the number of subdomains sharing $x$, and we simply take $\delta^{\dagger}(x) = 1/\N_x$. In applications, these scaling factors will depend on the heat conduction coefficient and the first of the Lam\'{e} parameters for scalar elliptic problems and the equations of linear
elasticity, respectively; see \cite{kla02,kla06}. Given such scaling factors at the subdomain boundary nodes, we can define a scaled operator $B_{\Delta, D}$. We note that each row of $B_\Delta$ has only two nonzero entries, $1$ and $-1$, corresponding to the same subdomain boundary node $x$. Multiplying each entry by the scaling factor $\delta^{\dagger}(x)$ gives us $B_{\Delta, D}$.

Solving the original fully assembled linear system~\EQ{matrix} is then equivalent to: find
$\left( {\bf u}_I, ~p_I, ~{\bf u}_{\Delta}, ~{\bf u}_{\Pi}, ~p_{\Gamma}, ~\lambda \right) \in
{\bf W}_I \bigoplus Q_I \bigoplus {\bf W}_{\Delta} \bigoplus {\bf W}_\Pi \bigoplus Q_\Gamma \bigoplus \Lambda$, such that
\be
\label{equation:bigeq}
\left[
\begin{array}{cccccc}
A_{II}      & B_{II}^T       & A_{I \Delta}       & A_{I \Pi}       & B_{\Gamma I}^T     &  0           \\[0.8ex]
B_{II}      & 0              & B_{I \Delta}       & B_{I \Pi}       & 0                  &  0           \\[0.8ex]
A_{\Delta I}& B_{I \Delta} ^T& A_{\Delta\Delta}   & A_{\Delta \Pi}  & B_{\Gamma \Delta}^T&  B_{\Delta}^T\\[0.8ex]
A_{\Pi I}   & B_{I \Pi}^T    & A_{\Pi \Delta}     & A_{\Pi \Pi}     & B_{\Gamma \Pi}^T   &  0           \\[0.8ex]
B_{\Gamma I}& 0              & B_{\Gamma \Delta}  & B_{\Gamma \Pi}  & 0                  &  0           \\[0.8ex]
0           & 0              & B_{\Delta}         & 0               & 0                  &  0
\end{array}
\right]
\left[ \begin{array}{c}
{\bf u}_I        \\[0.8ex]
p_I              \\[0.8ex]
{\bf u}_{\Delta} \\[0.8ex]
{\bf u}_{\Pi}    \\[0.8ex]
p_{\Gamma}     \\[0.8ex]
\lambda
\end{array} \right] =
\left[ \begin{array}{l}
{\bf f}_I        \\[0.8ex]
0                \\[0.8ex]
{\bf f}_{\Delta} \\[0.8ex]
{\bf f}_\Pi      \\[0.8ex]
0                \\[0.8ex]
0
\end{array} \right] \mbox{ ,  }
\ee
where the sub-blocks in the coefficient matrix represent the restrictions of $A$ and $B$ in~\EQ{matrix} to appropriate subspaces. The leading three-by-three block can be made block diagonal with each diagonal block representing one independent subdomain problem.

Corresponding to the one-dimensional null space of~\EQ{matrix}, we consider a vector of the form
$\left( \u_I,~p_I, ~\u_\Delta, ~ \u_\Pi, ~p_\Gamma, ~\lambda \right)  = \left( {\bf 0},~1_{p_I},~{\bf 0},
~{\bf 0}, ~1_{p_\Gamma}, \lambda \right)$, where $1_{p_I} \in Q_I$ and
$1_{p_\Gamma} \in Q_\Gamma$ represent vectors with value $1$ on each
entry. Substituting it into \EQ{bigeq} gives zero blocks on the  right-hand side, except at the third block
\be
\label{equation:fdelta}
{\bf f}_{\Delta} = [B_{I\Delta}^T ~~ B_{\Gamma\Delta}^T]\left[\begin{array}{c}1_{p_I}\\
1_{p_\Gamma}\end{array}\right]+B^T_\Delta \lambda.
\ee
The first term on the right-hand side represents the line integral of the normal component of the velocity finite element basis functions across the subdomain boundary on neighboring subdomains. Corresponding to the same subdomain boundary velocity degree of freedom, their values on the two neighboring subdomains are negative of each other. Therefore
\[
[B_{I\Delta}^T ~~ B_{\Gamma\Delta}^T]\left[\begin{array}{c}1_{p_I}\\
1_{p_\Gamma}\end{array}\right] = B^T_\Delta B_{\Delta,D}[B_{I\Delta}^T ~~ B_{\Gamma\Delta}^T]\left[\begin{array}{c}1_{p_I}\\
1_{p_\Gamma}\end{array}\right],
\]
from which we know that $\vvec{f}_\Delta = {\bf 0}$, for
\[
\lambda =-B_{\Delta,D}[B_{I\Delta}^T ~~ B_{\Gamma\Delta}^T]\left[\begin{array}{c}1_{p_I}\\
   1_{p_\Gamma}\end{array}\right].
\]
Therefore,  a basis of the one-dimensional null space of \EQ{bigeq} is
\begin{equation}
\label{equation:bignull}
\left( \begin{array}{cccccc}
0, & 1_{p_I}, & 0, & 0, & 1_{p_\Gamma}, &
-B_{\Delta,D}[B_{I\Delta}^T ~~ B_{\Gamma\Delta}^T]\left[\begin{array}{c}1_{p_I}\\ 1_{p_\Gamma}\end{array}\right] \end{array} \right).
\end{equation}

\section{A reduced symmetric positive semi-definite system}
\label{section:Gmatrix}

The system \EQ{bigeq} can be reduced to a Schur complement problem for the variables $\left(p_{\Gamma}, ~\lambda \right)$.
Since the leading four-by-four block of the coefficient matrix in \EQ{bigeq} is invertible, the variables  $\left( {\bf u}_I, ~p_I, ~{\bf u}_{\Delta},
~{\bf u}_{\Pi} \right)$ can be eliminated and we obtain
\begin{equation}
\label{equation:spd}
G \left[ \begin{array}{c}
p_\Gamma         \\[0.8ex]
\lambda
\end{array} \right] ~ = ~ g,
\end{equation}
where
\begin{equation}
\label{equation:Gmatrix}
G = \left[
\begin{array}{cccc}
B_{\Gamma I} & 0 & B_{\Gamma \Delta} & B_{\Gamma \Pi} \\[0.8ex]
0            & 0 & B_{\Delta}        & 0              \end{array}
\right] \left[
\begin{array}{cccc}
A_{II}       & B_{II}^T        & A_{I \Delta}      & A_{I \Pi}       \\[0.8ex]
B_{II}       & 0               & B_{I \Delta}      & B_{I \Pi}       \\[0.8ex]
A_{\Delta I} & B_{I \Delta} ^T & A_{\Delta\Delta}  & A_{\Delta \Pi}  \\[0.8ex]
A_{\Pi I}    & B_{I \Pi}^T     & A_{\Pi \Delta}    & A_{\Pi \Pi}
\end{array}
\right]^{-1} \left[
\begin{array}{cc}
B_{\Gamma I}^T        & 0            \\[0.8ex]
0                     & 0            \\[0.8ex]
B_{\Gamma \Delta}^T   & B_{\Delta}^T \\[0.8ex]
B_{\Gamma \Pi}^T      & 0
\end{array}
\right],
\end{equation}
and
\begin{equation}
\label{equation:gvec}
g = \left[
\begin{array}{cccc}
B_{\Gamma I} & 0 & B_{\Gamma \Delta} & B_{\Gamma \Pi} \\[0.8ex]
0            & 0 & B_{\Delta}        & 0
\end{array} \right] \left[
\begin{array}{cccc}
A_{II}       & B_{II}^T        & A_{I \Delta}      & A_{I \Pi}       \\[0.8ex]
B_{II}       & 0               & B_{I \Delta}      & B_{I \Pi}       \\[0.8ex]
A_{\Delta I} & B_{I \Delta} ^T & A_{\Delta\Delta}  & A_{\Delta \Pi}  \\[0.8ex]
A_{\Pi I}    & B_{I \Pi}^T     & A_{\Pi \Delta}    & A_{\Pi \Pi}
\end{array}
\right]^{-1} \left[
\begin{array}{l}
{\bf f}_I        \\[0.8ex]
0                \\[0.8ex]
{\bf f}_{\Delta} \\[0.8ex]
{\bf f}_\Pi
\end{array}
\right].
\end{equation}

We denote
\begin{equation}
\label{equation:AtildeBc}
\widetilde{A} = \left[
\begin{array}{cccc}
A_{II}      & B_{II}^T       & A_{I \Delta}       & A_{I \Pi}       \\[0.8ex]
B_{II}      & 0              & B_{I \Delta}       & B_{I \Pi}       \\[0.8ex]
A_{\Delta I}& B_{I \Delta} ^T& A_{\Delta\Delta}   & A_{\Delta \Pi}  \\[0.8ex]
A_{\Pi I}   & B_{I \Pi}^T    & A_{\Pi \Delta}     & A_{\Pi \Pi}   \end{array}
\right]
\quad \mbox{and}\quad
B_C=\left[
\begin{array}{cccc}
B_{\Gamma I} & 0 & B_{\Gamma \Delta} & B_{\Gamma \Pi} \\[0.8ex]
0            & 0 & B_{\Delta}        & 0              \end{array}
\right].
\end{equation}
We can see that $-G$ is the Schur complement of the coefficient matrix of \EQ{bigeq} with respect to the last two row blocks, i.e.,
\[
\left[ \begin{array}{cc} I & 0 \\[0.8ex] -B_C \widetilde{A}^{-1} & I \end{array} \right]
\left[ \begin{array}{cc} \widetilde{A} & B_C^T \\[0.8ex] B_C & 0 \end{array} \right]
\left[ \begin{array}{cc} I & - \widetilde{A}^{-1} B_C^T \\[0.8ex] 0  & I \end{array} \right] =
\left[ \begin{array}{cc} \widetilde{A} & 0 \\[0.8ex] 0 & -G \end{array} \right].
\]
From the Sylvester's law of inertia, namely, the number of positive, negative, and zero eigenvalues of a symmetric matrix is invariant under a change of coordinates, we can see that the number of zero eigenvalues of $G$ is the same as the number of zero eigenvalues (with multiplicity counted) of the original coefficient matrix of \EQ{bigeq}, which is one, and all other eigenvalues of $G$ are positive. Therefore $G$ is symmetric positive semi-definite. The null space of $G$ is derived from the null space of the original coefficient matrix of \EQ{bigeq}, and its basis is given by, cf.~\EQ{bignull},
\[
\left( \begin{array}{cc} 1_{p_\Gamma}, & - B_{\Delta,D}[B_{I\Delta}^T ~~ B_{\Gamma\Delta}^T]\left[\begin{array}{c}1_{p_I}\\ 1_{p_\Gamma} \end{array} \right]
\end{array} \right).
\]

We denote $X = Q_\Gamma \bigoplus \Lambda$. The range of $G$, denoted by $R_G$, is the subspace of $X$  orthogonal to the null space of $G$, and has the form
\be
\label{equation:Grange}
R_G=\left\{  \left[ \begin{array}{c}
g_{p_\Gamma}         \\[0.8ex]
g_{\lambda}
\end{array} \right] \in X ~ {\Big|} ~ g_{p_\Gamma}^T 1_{{p_\Gamma}} -
g_{\lambda}^T \left(B_{\Delta,D}[B_{I\Delta}^T ~~ B_{\Gamma\Delta}^T]\left[\begin{array}{c}1_{{p_I}}\\ 1_{{p_\Gamma}}\end{array}\right]\right)=0\right\}.
\ee

The restriction of $G$ to its range $R_G$ is positive definite.
The fact that the solution of \EQ{bigeq} always exists for any given $\left( {\bf f}_I, ~{\bf f}_{\Delta}, ~{\bf f}_{\Pi} \right)$ on the right-hand side implies that the solution of~\EQ{spd} exits for any $g$ defined by \EQ{gvec}. Therefore $g \in R_G$. When the conjugate gradient method (CG) is applied to solve \EQ{spd} with zero initial guess, all the iterates are in the Krylov subspace generated by $G$ and $g$, which is also a subspace of $R_G$, and where the CG cannot break down. After obtaining $\left( p_{\Gamma}, ~\lambda \right)$ from solving \EQ{spd}, the other components $\left( {\bf u}_I, ~p_I, ~{\bf u}_{\Delta}, ~{\bf u}_{\Pi} \right)$ in \EQ{bigeq} are obtained by back substitution.

In the rest of this section, we discuss the implementation of multiplying $G$ by a vector. The main operation is the product of $\Atilde^{-1}$ with a vector, cf. \EQ{Gmatrix} and \EQ{gvec}.  We denote
\[ A_{rr} = \left[ \begin{array}{ccc}
A_{II}       & B_{II}^T        & A_{I \Delta}     \\[0.8ex]
B_{II}       & 0               & B_{I \Delta}     \\[0.8ex]
A_{\Delta I} & B_{I \Delta} ^T & A_{\Delta\Delta} \end{array} \right] ,  \quad
A_{\Pi r} = A_{r \Pi}^T = \left[ A_{\Pi I}  \quad B_{I \Pi}^T  \quad A_{\Pi \Delta} \right], \quad  f_r = \left[ \begin{array}{l}
{\bf f}_I        \\[0.8ex]
0                \\[0.8ex]
{\bf f}_{\Delta} \end{array} \right],
\]
and define the Schur complement
\[
S_{\Pi} = A_{\Pi \Pi} - A_{\Pi r} A_{rr}^{-1} A_{r \Pi},
\]
which is symmetric positive definite from the Sylvester's law of inertia. $S_\Pi$ defines the coarse level problem in the algorithm.  The product
\[
   \left[
   \begin{array}{cccc}
   A_{II}       & B_{II}^T        & A_{I \Delta}      & A_{I \Pi}      \\[0.8ex]
   B_{II}       & 0               & B_{I \Delta}      & B_{I \Pi}      \\[0.8ex]
   A_{\Delta I} & B_{I \Delta} ^T & A_{\Delta\Delta}  & A_{\Delta \Pi} \\[0.8ex]
   A_{\Pi I}    & B_{I \Pi}^T     & A_{\Pi \Delta}    & A_{\Pi \Pi}
   \end{array}
   \right]^{-1} \left[
   \begin{array}{l}
   {\bf f}_I        \\[0.8ex]
   0                \\[0.8ex]
   {\bf f}_{\Delta} \\[0.8ex]
   {\bf f}_\Pi
   \end{array}
   \right]
\]
can then be represented by
\[
\left[ \begin{array}{c} A_{rr}^{-1} f_r \\[0.8ex] \vvec{0} \end{array} \right] ~ + ~
\left[ \begin{array}{c} -A_{rr}^{-1} A_{r \Pi} \\[0.8ex] I_\Pi \end{array}  \right] ~ S_{\Pi}^{-1} ~
\left({\bf f}_\Pi - A_{\Pi r} A_{rr}^{-1} f_r \right),
\]
which requires solving the coarse level problem once and independent subdomain Stokes problems with Neumann type boundary conditions twice.

\section{Some techniques}
\label{section:techniques}

We first define certain norms for several vector/function spaces. We denote
\be
\label{equation:Wtilde}
\vvec{\Wtilde} = {\bf W}_I \bigoplus \vvec{\Wtilde}_{\Gamma}.
\ee
For any $\vvec{w}$ in $\vvec{\Wtilde}$, we denote its restriction to subdomain $\Omega_i$ by ${\bf w}^{(i)}$. A subdomain-wise $H^1$-seminorm can be defined for functions in $\vvec{\Wtilde}$ by
\[
|\vvec{w}|^2_{H^1} = \sum_{i=1}^N |\vvec{w}^{(i)}|^2_{H^1(\Omega_i)}.
\]

We also define
\[
\Wtilde = {\bf W}_I \bigoplus Q_I \bigoplus {\bf W}_{\Delta} \bigoplus {\bf W}_\Pi,
\]
and its subspace
\be
\label{equation:W0}
\Wtilde_0 = \left\{ w = \left( {\bf w}_I, ~p_I, ~{\bf w}_{\Delta}, ~{\bf w}_{\Pi} \right) \in \Wtilde ~ \big| ~  B_{I I} \vvec{w}_I + B_{I\Delta} \vvec{w}_\Delta + B_{I\Pi} \vvec{w}_\Pi = 0 \right\}.
\ee
For any $w = \left( {\bf w}_I, ~p_I, ~{\bf w}_{\Delta}, ~{\bf w}_{\Pi} \right) \in \Wtilde_0$, let
$\vvec{w} =  \left( {\bf w}_I, ~{\bf w}_{\Delta}, ~{\bf w}_{\Pi} \right) \in \vvec{\Wtilde}$. Then
\begin{eqnarray}
\left< w, w \right>_{\widetilde{A}} & = &
\left[ \begin{array}{l} {\bf w}_I \\[0.8ex] {\bf w}_{\Delta} \\[0.8ex] {\bf w}_\Pi \end{array} \right]^T
\left[ \begin{array}{ccc}
A_{II}        & A_{I \Delta}      & A_{I \Pi}      \\[0.8ex]
A_{\Delta I}  & A_{\Delta\Delta}  & A_{\Delta \Pi} \\[0.8ex]
A_{\Pi I}     & A_{\Pi \Delta}    & A_{\Pi \Pi} \end{array} \right]
\left[ \begin{array}{l} {\bf w}_I \\[0.8ex] {\bf w}_{\Delta} \\[0.8ex] {\bf w}_\Pi \end{array} \right] \nonumber \\[0.8ex]
& = & \sum_{i=1}^N \left[ \begin{array}{c} {\bf w}_I^{(i)} \\ {\bf w}_{\Delta}^{(i)} \\ {\bf w}_{\Pi}^{(i)} \end{array} \right]^T
\left[ \begin{array}{cccc}
A_{II}^{(i)}       & A_{I \Delta}^{(i)}      & A_{I \Pi}^{(i)}      \\[0.8ex]
A_{\Delta I}^{(i)} & A_{\Delta\Delta}^{(i)}  & A_{\Delta \Pi}^{(i)} \\[0.8ex]
A_{\Pi I}^{(i)}    & A_{\Pi \Delta}^{(i)}    & A_{\Pi \Pi}^{(i)}
\end{array} \right] \left[ \begin{array}{c} {\bf w}_I^{(i)} \\ {\bf w}_{\Delta}^{(i)} \\ {\bf w}_{\Pi}^{(i)} \end{array} \right]
= \sum_{i=1}^N \left| \left[ \begin{array}{c} {\bf w}_I^{(i)} \\ {\bf w}_{\Delta}^{(i)} \\ {\bf w}_{\Pi}^{(i)} \end{array} \right] \right|_{H^1(\Omega_i)}^2 \label{equation:W0n} \\[0.8ex]
& = & |\vvec{w}|^2_{H^1},  \nonumber
\end{eqnarray}
i.e., $\left< \cdot , \cdot \right>_{\Atilde}$ defines an inner product on $\Wtilde_0$. In \EQ{W0n}, the superscript ${}^{(i)}$ is used to represent the restrictions of corresponding vectors and matrices to subdomain $\Omega_i$.

Since $\vvec{W}$ is essentially the subspace of $\vvec{\Wtilde}$ with continuous subdomain boundary velocities, the inf-sup condition \EQ{infsup} and \EQ{infsupMatrix} also holds for the mixed space $\vvec{\Wtilde} \times Q$. Denote
\begin{equation}\label{equation:Btilde}
\widetilde{B} = \left[ \begin{array}{ccc}
B_{II}      & B_{I \Delta}       & B_{I \Pi}      \\[0.8ex]
B_{\Gamma I}& B_{\Gamma \Delta}  & B_{\Gamma \Pi}
\end{array} \right], \qquad
\overline{\widetilde{A}} = \left[
\begin{array}{ccc}
A_{II}      & A_{I \Delta}       & A_{I \Pi}       \\[0.8ex]
A_{\Delta I}& A_{\Delta\Delta}   & A_{\Delta \Pi}  \\[0.8ex]
A_{\Pi I}   & A_{\Pi \Delta}     & A_{\Pi \Pi}   \end{array}
\right],
\end{equation}
as in \EQ{bigeq}, then
\begin{equation}
\label{equation:infsupMatrixtilde} \sup_{{\bf w} \in \vvec{\Wtilde}}
\frac{\left< q, \widetilde{B} \vvec{w} \right>^2}{\left< \vvec{w},  \overline{\widetilde{A}} \vvec{w} \right>} \geq \beta^2 \left< q, Z q \right>,
\hspace{0.5cm} \forall q \in Q/Ker(B^T),
\end{equation}
where $\beta$ is the same as in \EQ{infsup} and \EQ{infsupMatrix}.

We also have the following lemma on the stability of the operator $\widetilde{B}$.

\begin{mylemma}
\label{lemma:BtildeStability}
For any $\vvec{w} \in \vvec{\Wtilde}$ and $q \in Q$, $\left< \Btilde {\bf w}, q \right> \leq | \vvec{w} |_{H^1} \| q \|_{L^2}$.
\end{mylemma}

\beginproof
\begin{eqnarray*}
\left<\Btilde {\bf w}, q \right>^2 & = & \left( \sum_{i=1}^N \int_{\Omega_i} \nabla\cdot {\bf w}^{(i)} q \right)^2\leq \left( \sum_{i=1}^N \sqrt{\int_{\Omega_i} | \nabla {\bf w}^{(i)} |^2} \sqrt{\int_{\Omega_i} q^2} \right)^2  \\[0.8ex]
& \le & \left( \sum_{i=1}^N \int_{\Omega_i} | \nabla {\bf w}^{(i)} |^2 \right) \left( \sum_{i=1}^N \int_{\Omega_i} q^2 \right) = | \vvec{w} |^2_{H^1} \| q \|^2_{L^2}. \qquad \Box
\end{eqnarray*}

The finite element space for subdomain boundary pressures, $Q_\Gamma$, is a subspace of $L^2(\Gamma)$. For each $p_\Gamma \in Q_\Gamma$, its finite element extension by zero to the interior of subdomains is denoted by $p_\Gamma^E$, which equals $p_\Gamma$ on all subdomain boundary nodes and equals zero on all subdomain interior nodes. We can see that $p_\Gamma^E \in Q \subset L^2(\Omega)$, and $\| p_\Gamma^E \|_{L^2(\Omega)}^2 = \left< p^E_\Gamma, p^E_\Gamma \right>_Z$, from the definition of $Z$ in Section~\ref{section:FEM}.

From \EQ{Gmatrix} and \EQ{AtildeBc}, we can see that
\[
G = B_C \Atilde^{-1} B_C^T.
\]
In particular, we denote the first row of $B_C$ by
\[
\widetilde{B}_{\Gamma} =
\left[ B_{\Gamma I} \quad 0  \quad B_{\Gamma \Delta}  \quad  B_{\Gamma \Pi} \right];
\]
for the second row, we denote the restriction operator from $\Wtilde$ onto ${\bf W}_{\Delta}$ by $\widetilde{R}_\Delta$, such that for any $w = \left( {\bf w}_I, ~p_I, ~{\bf w}_{\Delta}, ~{\bf w}_{\Pi} \right) \in \Wtilde$,
$\widetilde{R}_\Delta w = {\bf w}_{\Delta}$.
Then $G$ can be represented by the following two-by-two block structure
\begin{equation}
\label{equation:Gtwo}
G = \left[ \begin{array}{cc} G_{p_\Gamma p_\Gamma} & G_{p_\Gamma \lambda} \\[0.8ex] G_{\lambda p_\Gamma} & G_{\lambda \lambda} \end{array} \right],
\end{equation}
where
\begin{eqnarray*}
& G_{p_\Gamma p_\Gamma} = \widetilde{B}_{\Gamma} \widetilde{A}^{-1} \widetilde{B}_{\Gamma}^T, \qquad
G_{p_\Gamma \lambda} = \widetilde{B}_{\Gamma} \widetilde{A}^{-1} \widetilde{R}_{\Delta}^T B_{\Delta}^T, & \\ [0.8ex]
& G_{\lambda p_\Gamma} = B_{\Delta} \widetilde{R}_\Delta \widetilde{A}^{-1} \widetilde{B}_{\Gamma}^T, \qquad
G_{\lambda \lambda} = B_{\Delta} \widetilde{R}_\Delta \widetilde{A}^{-1} \widetilde{R}_{\Delta}^T B_{\Delta}^T. &
\end{eqnarray*}

The pressure components of all  vectors in $R_G$ with $g_\lambda = 0$, cf. \EQ{Grange}, form a subspace of $Q_\Gamma$ and we denote this subspace by $R_{G|Q_\Gamma}$. From the definition of $R_G$, we can see that for any vector $p_\Gamma \in R_{G|Q_\Gamma}$, $p_\Gamma^T 1_{p_\Gamma} = 0$, and then its extension by zero to the interior of subdomains, $p_\Gamma^E$, also has zero average.

The following lemma follows essentially from \cite[Lemma 9.1]{Toselli:2004:DDM}.
\begin{mylemma}
\label{lemma:mass}
For all $p_\Gamma \in R_{G|Q_\Gamma}$,
\[
\beta^2 \| p_\Gamma^E \|_{L^2(\Omega)}^2 ~ \leq ~ \left<p_\Gamma, G_{p_\Gamma p_\Gamma} p_\Gamma \right> ~ \leq ~
\| p_\Gamma^E \|_{L^2(\Omega)}^2,
\]
where $p_\Gamma^E$ represents the extension by zero of $p_\Gamma$ to the interior of subdomains, and $\beta$ is the same as in \EQ{infsup} and \EQ{infsupMatrix}.
\end{mylemma}

\beginproof
Note that even though $\widetilde{A}^{-1}$ is indefinite in $\Wtilde$, it is positive definite when restricted to a subspace of $\Wtilde$, where the pressure component equals zero, and the norm $\| \cdot \|_{\widetilde{A}^{-1}}$ is well defined.

To prove the left side inequality, denote for any ${\bf v} = \left( {\bf v}_I,  ~{\bf v}_{\Delta}, ~{\bf v}_{\Pi} \right) \in \vvec{\Wtilde}$, ${\bf v}^\dagger = \left( {\bf v}_I,  ~ 0, ~{\bf v}_{\Delta}, ~{\bf v}_{\Pi} \right) \in \Wtilde$. We have
\begin{eqnarray*}
\left<p_\Gamma, \widetilde{B}_{\Gamma} \widetilde{A}^{-1} \widetilde{B}_{\Gamma}^T p_\Gamma \right> = \| \widetilde{B}_{\Gamma}^T p_\Gamma \|_{\widetilde{A}^{-1}}^2 = \sup_{{\bf v} \in \vvec{\Wtilde}}
\frac{\left< {\bf v}^\dagger, \widetilde{B}_{\Gamma}^T p_\Gamma \right>^2_{\widetilde{A}^{-1}}}{\| {\bf v}^\dagger \|^2_{\widetilde{A}^{-1}}}
= \sup_{{\bf v} \in \vvec{\Wtilde}}
\frac{\left( p_\Gamma^T \widetilde{B}_{\Gamma} \widetilde{A}^{-1} {\bf v}^\dagger \right)^2}{{\bf v}^{\dagger^T} {\widetilde{A}^{-1}} {\bf v}^\dagger}  \\
= \sup_{{\bf w} \in \vvec{\Wtilde}}
\frac{\left( p_\Gamma^{T} \widetilde{B}_\Gamma {\bf w}^\dagger \right)^2}{{\bf w}^{\dagger^T} \widetilde{A} {\bf w}^\dagger} = \sup_{{\bf w} \in \vvec{\Wtilde}}
\frac{\left( p_\Gamma^{E^T} \widetilde{B} {\bf w} \right)^2}{{\bf w}^T \overline{\widetilde{A}} {\bf w}} \geq \beta^2 \left< p_\Gamma^E, p_\Gamma^E \right>_Z = \beta^2 \| p_\Gamma^E \|_{L^2(\Omega)}^2 ,
\end{eqnarray*}
where we have used the inf-sup condition \EQ{infsupMatrixtilde} for the inequality in the middle.

To prove the right side inequality, for any given $p_\Gamma \in R_{G|Q_\Gamma}$, denote ${\bf v}^\dagger = \left( {\bf v}_I,  ~ p_I, ~{\bf v}_{\Delta}, ~{\bf v}_{\Pi} \right) = \widetilde{A}^{-1} \widetilde{B}_{\Gamma}^T p_\Gamma$, and the shorter vector ${\bf v} = \left( {\bf v}_I, ~{\bf v}_{\Delta}, ~{\bf v}_{\Pi} \right)$. From the continuity of $\Btilde$ in Lemma \LA{BtildeStability} and \EQ{W0n}, we have
\begin{eqnarray*}
& & \left<p_\Gamma, \widetilde{B}_{\Gamma} \widetilde{A}^{-1} \widetilde{B}_{\Gamma}^T p_\Gamma \right>  = \left<p_\Gamma, \widetilde{B}_{\Gamma} {\bf v}^\dagger \right> = \left<p_\Gamma^E, \widetilde{B} \vvec{v} \right> \leq \| p_\Gamma^E \|_{L^2} ~ | \vvec{v} |_{H^1} \\
& = & \| p_\Gamma^E \|_{L^2} ~ \sqrt{\left< \widetilde{A}^{-1} \widetilde{B}_{\Gamma}^T p_\Gamma, \widetilde{A}^{-1} \widetilde{B}_{\Gamma}^T p_\Gamma \right>_{\Atilde}} = \| p_\Gamma^E \|_{L^2}  ~ \left<p_\Gamma, \widetilde{B}_{\Gamma} \widetilde{A}^{-1} \widetilde{B}_{\Gamma}^T p_\Gamma \right>^{1/2}. \qquad \Box
\end{eqnarray*}

The following corollary of Lemma~\ref{lemma:mass} is an immediate result from \EQ{massmatrix} and the facts that $\| p_\Gamma^E \|_{L^2(\Omega)}^2 = \left< p_\Gamma^E, p_\Gamma^E \right>_Z$, $\left<p_\Gamma^E, p_\Gamma^E \right> = \left<p_\Gamma, p_\Gamma \right>$.
\begin{mycorollary}
\label{coro:mass}
There exist positive constants $c$ and $C$, such that
\[
c h^2 \beta^2 I_{p_\Gamma} ~ \leq ~ G_{p_\Gamma p_\Gamma}  ~ \leq ~ C h^2  I_{p_\Gamma}
\]
where $I_{p_\Gamma}$ is the identity matrix of the same dimension as $G_{p_\Gamma p_\Gamma}$, and $\beta$ is the same as in \EQ{infsup} and \EQ{infsupMatrix}.
\end{mycorollary}

\begin{myremark}
Lemma \ref{lemma:mass} and Corollary \ref{coro:mass} are not used in our proof of the condition number bound in Section~\ref{section:lumped}. However, it is intuitive to see from Corollary~\ref{coro:mass} that the first diagonal block $G_{p_\Gamma p_\Gamma}$ in matrix $G$ can be approximated spectrally equivalently by the identity matrix multiplied by $h^2$, which is what is being done in our block diagonal preconditioner discussed in Section~\ref{section:lumped}.
\end{myremark}

We also need define a certain jump operator across the subdomain boundaries $\Gamma$. Let $P_D: \Wtilde \rightarrow \Wtilde$, be defined by, cf.~\cite{li06A},
\[
P_D = \widetilde{R}_{\Delta}^T B_{\Delta, D}^T B_{\Delta} \widetilde{R}_\Delta.
\]
We can see that application of $P_D$ to a vector essentially computes the difference (jump) of the dual velocity component across the subdomain boundaries and then distributes the jump to neighboring subdomains according to the scaling factor $\delta^\dagger(x)$. In fact, the dual velocity component is the only part of the vector involved in the application of $P_D$; all other components are kept zero and are added into the definition to make $P_D$ more convenient to use in the presentation of the algorithm. We also have, for any $w = \left( {\bf w}_I, ~p_I, ~{\bf w}_{\Delta}, ~{\bf w}_{\Pi} \right) \in \Wtilde$,
\[
\left< P_D w, P_D w \right>_{\widetilde{A}} = \left< B_{\Delta, D}^T B_{\Delta} {\bf w}_{\Delta}, B_{\Delta, D}^T B_{\Delta} {\bf w}_{\Delta} \right>_{A_{\Delta \Delta}}.
\]
The following lemma can be found essentially from \cite[Section 6]{li07}; see also \EQ{W0n}.
\begin{mylemma}
\label{lemma:jump}
There exists a function $\Phi(H/h)$, such that for all $w \in \Wtilde_0$,
\[
\left< P_D w, P_D w \right>_{\widetilde{A}} \leq  \Phi(H/h) \left< w, w \right>_{\widetilde{A}}.
\]
\end{mylemma}

\begin{myremark}
\label{remark:Phi}
Just as for the positive definite elliptic problems discussed in \cite[Section~6]{li07}, for two-dimensional problems, when only subdomain corner velocities are chosen as coarse level primal variables, $\Phi(H/h) = C (H/h) (1 + \log{(H/h)})$; when both subdomain corner and edge-average velocity degrees of freedom are chosen as primal variables, $\Phi(H/h) = C H/h$.
\end{myremark}

The following lemma is also used and can be found at \cite[Lemma~2.3]{paulo2003}.
\begin{mylemma}
\label{lemma:paul}
Consider the saddle point problem: find  $(\vvec{u}, p) \in \vvec{W} \bigoplus Q$, such that
\be \left[
\begin{array}{cc}
A       & B^T      \\[0.8ex]
B       & 0
\end{array}
\right] \left[
\begin{array}{l}
\vvec{u}        \\[0.8ex]
p
\end{array}
\right]
=\left[
\begin{array}{l}
\vvec{f}        \\[0.8ex]
g
\end{array}
\right],
\ee
where $A$ and $B$ are as in \EQ{matrix}, $\vvec{f} \in \vvec{W}$, and $g \in Im(B) \subset Q$. Let $\beta$ be the inf-sup constant specified in \EQ{infsupMatrix}.
Then
\[
\| \vvec{u} \|_A \le \| \vvec{f} \|_{A^{-1}} + \frac{1}{\beta} \| g \|_{Z^{-1}},
\]
where $Z$ is the mass matrix defined in Section \ref{section:FEM}.
\end{mylemma}

\section{A lumped preconditioner}
\label{section:lumped}

The lumped preconditioner was first used in the FETI algorithm
\cite{FETI2} for solving positive definite elliptic problems. Compared
with the Dirichlet preconditioner, also used for the FETI algorithm
\cite{FETI4}, the lumped preconditioner is less effective in the
improvement of convergence rate, but it is also less expensive in the
computational costs. The main operation in the lumped preconditioner
is subdomain matrix and vector products, while the implementation of
the Dirichlet preconditioner requires solving subdomain
systems of equations. In this paper, we discuss only the lumped
preconditioner in our algorithm for solving the incompressible Stokes
equation; study of the Dirichlet preconditioner will be addressed in
forthcoming work.

We consider a block diagonal preconditioner for \EQ{spd}. From Corollary~\ref{coro:mass}, the inverse of the first diagonal block $G_{p_\Gamma p_\Gamma}$ of $G$ can be effectively approximated by $1/h^2$ times the identity matrix. The inverse of the second diagonal block $B_{\Delta} \widetilde{R}_\Delta \widetilde{A}^{-1} \widetilde{R}_{\Delta}^T B_{\Delta}^T$, can be approximated by the following lumped block
\[
M^{-1}_\lambda = B_{\Delta, D} \widetilde{R}_\Delta \widetilde{A} \widetilde{R}_{\Delta}^T B_{\Delta, D}^T.
\]
This leads to the lumped preconditioner
\[
M^{-1} = \left[ \begin{array}{cc} \frac{1}{h^2} I_{p_\Gamma} & \\[0.8ex]
& M^{-1}_\lambda  \end{array} \right],
\]
for solving \EQ{spd}.

\begin{myremark}
The mesh size $h$ is used in the above preconditioner. For applications where the mesh size is not explicitly provided and only the coefficient matrix in~\EQ{matrix} is given, an estimate of $h$ can be obtained by comparing the nonzero entries in $A$ and $B$ blocks. From the definition of $A$ and $B$ for the incompressible Stokes problem \EQ{bilinear}, entries in $A$ and entries in $B$ have a difference of factor $h$ in general.
\end{myremark}

$M^{-1}$ is symmetric positive definite. Multiplication of $M^{-1}$ by a vector requires mainly the product of $\widetilde{A}$ with a vector. When the CG iteration is applied to solve the preconditioned system
\begin{equation}
\label{equation:Mspd}
M^{-1} G \left[ \begin{array}{c}
p_\Gamma         \\[0.8ex]
\lambda
\end{array} \right] ~ = ~ M^{-1} g,
\end{equation}
with zero initial guess, all the iterates belong to the Krylov
subspace generated by the operator $M^{-1} G$ and the vector $M^{-1}
g$, which is also a  subspace of the range of $M^{-1} G$. We denote the range of $M^{-1} G$ by $R_{M^{-1} G}$. The following lemma shows that the CG iteration applied to solving \EQ{Mspd} cannot break down.

\begin{mylemma}
\label{lemma:CG} Let the preconditioner $M^{-1}$ be symmetric positive definite.
The CG iteration applied to solving \EQ{Mspd} with zero initial guess cannot break down.
\end{mylemma}

\beginproof
We just need to show that for any $0 \neq x \in R_{M^{-1} G}$, $G x \neq 0$. Let $0 \neq x = M^{-1} G y$, for a certain $y \in X$ and $y \neq 0$. $G x  = G M^{-1} G y$,
which cannot be zero since $G y \neq 0$ and $y^T  G M^{-1} G y \neq 0$. $\qquad \Box$

\begin{mylemma}
\label{lemma:m1RG} Let $M^{-1}$ be symmetric positive definite.
For any $x = (p_{\Gamma}, ~\lambda) \in R_{M^{-1} G}$,
\[
\left<Mx,x \right> = \max_{y \in R_G, y \neq 0} \frac{\left<y,x\right>^2}{\left<M^{-1}y,y\right>}.
\]
\end{mylemma}

\beginproof Denote the range of $M^{-\frac12}G$ by $R_{M^{-1/2} G}$. For any $x \in R_{M^{-1} G}$,
\begin{eqnarray*}
\left<Mx,x \right> & = & \left<M^{\frac12} x, M^{\frac12} x \right> =
\max_{z \in R_{M^{-1/2} G}, z \neq 0} \frac{\left<M^{\frac12} x, z\right>^2}{\left<z, z\right>} \\
& = & \max_{y \in R_G, y \neq 0}
\frac{\left<M^{\frac12} x, M^{-\frac12} y\right>^2}{\left<M^{-\frac12} y, M^{-\frac12} y\right>}
= \max_{y \in R_G, y \neq 0} \frac{\left<y,x\right>^2}{\left<M^{-1}y,y\right>} ~ .  \qquad \Box
\end{eqnarray*}

In the following, we establish a condition number bound of the preconditioned operator $M^{-1} G$. We first have the following lemma.

\begin{mylemma}\label{lemma:upper}
For any $w\in \Wtilde_0$,
\[
\left<M^{-1}B_Cw,B_Cw\right>\le \Phi(H/h)\left<\Atilde w, w\right>,
\]
where $\Phi(H/h)$ is as defined in Lemma \LA{jump}.
\end{mylemma}

\beginproof
Given $w = \left( {\bf w}_I, ~q_I, ~{\bf w}_{\Delta}, ~{\bf w}_{\Pi} \right) \in \Wtilde_0$, let
$g_{p_\Gamma} = B_{\Gamma I} {\bf w}_I + B_{\Gamma \Delta} {\bf w}_\Delta + B_{\Gamma\Pi} {\bf w}_\Pi$. We have
\BA
\label{equation:MBW}
\left<M^{-1}B_Cw,B_Cw\right> & = & \frac{1}{h^2}\left<g_{p_\Gamma}, g_{p_\Gamma} \right>+
\left(B_\Delta\Rtilde_\Delta w\right)^T M^{-1}_\lambda B_\Delta\Rtilde_\Delta w\nonumber\\
&=&\frac{1}{h^2}\left<g_{p_\Gamma}, g_{p_\Gamma} \right>+\left(B_\Delta\Rtilde_\Delta w\right)^T
B_{\Delta,D}\Rtilde_{\Delta}\Atilde\Rtilde_{\Delta}^TB_{\Delta,D}^T
\left(B_\Delta\Rtilde_\Delta w\right)\nonumber\\
&=&\frac{1}{h^2}\left< g_{p_\Gamma}, g_{p_\Gamma} \right>+\left<P_D w,P_D w\right>_\Atilde\nonumber\\
&\le&\frac{1}{h^2}\left< g_{p_\Gamma}, g_{p_\Gamma} \right>+\Phi(H/h)\left<w,w\right>_\Atilde,
\EA
where we used Lemma \LA{jump} for the last inequality. It is  sufficient to bound the first term of the right-hand side in the above inequality.

We denote $\vvec{w} = \left( {\bf w}_I, ~{\bf w}_{\Delta}, ~{\bf w}_{\Pi} \right) \in \vvec{\Wtilde}$. Since $B_{I I} \vvec{w}_I + B_{I\Delta} \vvec{w}_\Delta + B_{I\Pi} \vvec{w}_\Pi = 0$, cf.~\EQ{W0}, we have
\begin{eqnarray*}
\left<g_{p_\Gamma}, g_{p_\Gamma}\right> & = & \left[ \begin{array}{c} B_{I I} \vvec{w}_I + B_{I\Delta} \vvec{w}_\Delta + B_{I\Pi} \vvec{w}_\Pi \\ B_{\Gamma I} {\bf w}_I + B_{\Gamma \Delta} {\bf w}_\Delta + B_{\Gamma\Pi} {\bf w}_\Pi \end{array} \right]^T
\left[ \begin{array}{c} B_{I I} \vvec{w}_I + B_{I\Delta} \vvec{w}_\Delta + B_{I\Pi} \vvec{w}_\Pi \\ B_{\Gamma I} {\bf w}_I + B_{\Gamma \Delta} {\bf w}_\Delta + B_{\Gamma\Pi} {\bf w}_\Pi \end{array} \right] \\
& = & \left<\Btilde {\bf w},  \Btilde {\bf w} \right>,
\end{eqnarray*}
where  $\Btilde$ is defined in \EQ{Btilde}. From \EQ{massmatrix} and the stability of $\Btilde$, cf. Lemma \LA{BtildeStability}, we have
\begin{eqnarray}
\frac{1}{h^2}\left<g_{p_\Gamma}, g_{p_\Gamma}\right> & = & \frac{1}{h^2} \left<\Btilde {\bf w},  \Btilde {\bf w} \right> \leq C \left<\Btilde {\bf w},  \Btilde {\bf w} \right>_{Z^{-1}} =  C \max_{q \in Q} \frac{\left<\Btilde {\bf w}, q \right>^2}{\left<q, q\right>_Z} \label{equation:boundBu}\\
& \le & C \max_{q \in Q} \frac{| \vvec{w} |^2_{H^1}  \| q \|^2_{L^2}}{\| q \|^2_{L^2}} = C | \vvec{w} |^2_{H^1} = C \left<w,w\right>_\Atilde, \nonumber
\end{eqnarray}
where for the last equality, we used the fact that $B_{I I} \vvec{w}_I
+ B_{I\Delta} \vvec{w}_\Delta + B_{I\Pi} \vvec{w}_\Pi = 0$, and~\EQ{W0n}. $\quad \Box$

\begin{mylemma}
\label{lemma:lower}
For any given $y = (g_{p_{\Gamma}}, g_\lambda) \in R_G$, there exits
$w \in \Wtilde_0$, such that $B_C w = y$, and $\left <\Atilde w, w\right> \le \frac{C}{\beta^2}  \left< M^{-1}y, y \right>$.
\end{mylemma}

\beginproof
Given $y = (g_{p_{\Gamma}}, g_\lambda) \in R_G$, take ${\bf w}_{\Delta}^{(I)} = B_{\Delta, D}^T g_\lambda$.  Let $\w^{(I)} = ({\bf 0}, ~{\bf w}_{\Delta}^{(I)}, {\bf 0}) \in {\bf W}_I \bigoplus {\bf W}_{\Delta} \bigoplus {\bf W}_\Pi$ and $w^{(I)} = ( {\bf 0}, ~0, ~{\bf w}_{\Delta}^{(I)}, ~ {\bf 0} ) \in {\bf W}_I \bigoplus Q_I \bigoplus {\bf W}_{\Delta} \bigoplus {\bf W}_\Pi$. We have
\be
\label{equation:uOne}
| \w^{(I)} |^2_{H^1} =  \left< A_{\Delta\Delta}\w^{(I)}_\Delta,\w^{(I)}_\Delta\right>,
\ee
and
\be
\label{equation:bcWone}
B_c w^{(I)} = \left[ \begin{array}{cccc}
B_{\Gamma I} & 0 & B_{\Gamma \Delta} & B_{\Gamma \Pi} \\[0.8ex]
0            & 0 & B_{\Delta}        & 0              \end{array}
\right] \left[ \begin{array}{c}
{\bf 0} \\[0.8ex] 0 \\[0.8ex] B_{\Delta, D}^T g_\lambda \\[0.8ex] {\bf 0}
\end{array} \right] = \left[ \begin{array}{c}
B_{\Gamma\Delta}\w^{(I)}_\Delta \\[0.8ex] g_\lambda \end{array} \right],
\ee
where we used the fact that $B_\Delta B_{\Delta, D}^T = I$.

We consider the solution to the following fully assembled system of linear equations of the form~\EQ{matrix}: find
$({\bf w}_I^{(II)}, ~q_I^{(II)}, ~{\bf w}_{\Gamma}^{(II)}, ~q_\Gamma^{(II)}) \in {\bf W}_I \bigoplus Q_I \bigoplus {\bf W}_{\Gamma} \bigoplus Q_\Gamma$, such that
\be
\label{equation:uTwo}
\left[
\begin{array}{cccc}
A_{II}      & B_{II}^T       & A_{I \Gamma}      & B_{\Gamma I}^T     \\[0.8ex]
B_{II}      & 0              & B_{I \Gamma}      & 0                  \\[0.8ex]
A_{\Gamma I}& B_{I \Gamma} ^T& A_{\Gamma\Gamma}  & B_{\Gamma \Gamma}^T \\[0.8ex]
B_{\Gamma I}& 0              & B_{\Gamma \Gamma} & 0
\end{array}
\right]
\left[ \begin{array}{c}
{\bf w}_I^{(II)}        \\[0.8ex]
q_I^{(II)}              \\[0.8ex]
{\bf w}_{\Gamma}^{(II)} \\[0.8ex]
q_{\Gamma}^{(II)}
\end{array} \right] =
\left[ \begin{array}{l}
{\bf 0}        \\[0.8ex]
-B_{I\Delta}{\bf w}^{(I)}_\Delta              \\[0.8ex]
{\bf 0}        \\[0.8ex]
g_{p_{\Gamma}}-B_{\Gamma\Delta}{\bf w}^{(I)}_\Delta
\end{array} \right] \mbox{ , }
\ee
where a particular right-hand side is chosen. We first note that, since $(g_{p_{\Gamma}}, g_\lambda) \in R_G$, the right-hand side vector of the above system satisfies, cf. \EQ{Grange},
\[
( -B_{I\Delta}{\bf w}^{(I)}_\Delta )^T 1_{p_I}
+ ( g_{p_{\Gamma}}-B_{\Gamma\Delta}{\bf w}^{(I)}_\Delta )^T  1_{p_\Gamma}  =  g_{p_{\Gamma}}^T 1_{p_\Gamma} - g_\lambda^T B_{\Delta, D} \left( B_{I\Delta}^T 1_{p_I} + B_{\Gamma\Delta}^T 1_{p_\Gamma} \right) = 0,
\]
i.e., it has zero average, which implies existence of the solution to \EQ{uTwo}.

Denote ${\bf w}^{(II)} = ( {\bf w}_I^{(II)}, ~{\bf w}_{\Gamma}^{(II)} ) \in {\bf W}$. From the inf-sup stability of the original problem \EQ{matrix} and Lemma \ref{lemma:paul}, we have
\be \label{equation:uIbound}
| {\bf w}^{(II)} |^2_{H^1} \leq \frac{1}{\beta^2} \left\| \left[ \begin{array}{l}
-B_{I\Delta}\w^{(I)}_\Delta     \\[0.8ex] g_{p_\Gamma}-B_{\Gamma\Delta}\w^{(I)}_\Delta \end{array} \right] \right\|^2_{Z^{-1}} \le \frac{1}{\beta^2} \left\| \left[ \begin{array}{l} B_{I\Delta}\w^{(I)}_\Delta \\[0.8ex] B_{\Gamma\Delta}\w^{(I)}_\Delta \end{array} \right] \right\|^2_{Z^{-1}} + \frac{1}{\beta^2} \left\|\left[ \begin{array}{l} 0 \\[0.8ex] g_{p_\Gamma} \end{array} \right] \right\|^2_{Z^{-1}}.
\ee

The first term on the right-hand side of \EQ{uIbound} can be bounded in the same way as done in \EQ{boundBu}, and we have
\be
\label{equation:uuOne}
\left\| \left[ \begin{array}{l} B_{I\Delta}\w^{(I)}_\Delta \\[0.8ex] B_{\Gamma\Delta}\w^{(I)}_\Delta \end{array} \right] \right\|^2_{Z^{-1}} \leq C \left< A_{\Delta\Delta}\w^{(I)}_\Delta,\w^{(I)}_\Delta\right>;
\ee
the second term can be bounded by, using \EQ{massmatrix},
\be
\label{equation:uuTwo}
\left\| \left[ \begin{array}{l} 0   \\[0.8ex] g_{p_\Gamma} \end{array} \right] \right\|^2_{Z^{-1}}
\le \frac{C}{h^2} \left<g_{p_\Gamma},g_{p_\Gamma}\right>.
\ee

Split the continuous subdomain boundary velocity ${\bf w}_{\Gamma}^{(II)}$ into the dual part ${\bf w}_{\Delta}^{(II)} \in {\bf W}_{\Delta}$ and the primal part ${\bf w}_{\Pi}^{(II)} \in {\bf W}_{\Pi}$, and denote $w^{(II)} = ({\bf w}_I^{(II)}, ~q_I^{(II)}, ~{\bf w}_{\Delta}^{(II)}, ~{\bf w}_{\Pi}^{(II)})$. We have, from \EQ{uTwo},
\be
\label{equation:biWtwo}
\left[ \begin{array}{cccc} B_{II} & 0 & B_{I \Delta} & B_{I \Pi} \end{array} \right]
\left[ \begin{array}{c}
{\bf w}_I^{(II)}        \\[0.8ex]
q_I^{(II)}              \\[0.8ex]
{\bf w}_{\Delta}^{(II)} \\[0.8ex]
{\bf w}_{\Pi}^{(II)}
\end{array} \right] = -B_{I\Delta}\w^{(I)}_\Delta,
\ee
and
\be
\label{equation:bcWtwo}
B_c w^{(II)} = \left[
\begin{array}{cccc}
B_{\Gamma I} & 0 & B_{\Gamma \Delta} & B_{\Gamma \Pi} \\[0.8ex]
0            & 0 & B_{\Delta}        & 0              \end{array}
\right] \left[ \begin{array}{c}
{\bf w}_I^{(II)}        \\[0.8ex]
q_I^{(II)}              \\[0.8ex]
{\bf w}_{\Delta}^{(II)} \\[0.8ex]
{\bf w}_{\Pi}^{(II)}
\end{array} \right] = \left[ \begin{array}{c}
g_{p_\Gamma}    -B_{\Gamma\Delta}\w^{(I)}_\Delta   \\[0.8ex] 0  \end{array} \right].
\ee

Let
$w = w^{(I)} + w^{(II)}$. We can see from \EQ{biWtwo} that $w \in \Wtilde_0$, cf.~\EQ{W0}.
We can also see from \EQ{bcWone} and \EQ{bcWtwo} that $B_C w=
y$. Furthermore, by \EQ{W0n},
\[
| w |^2_{\widetilde{A}}=| \w^{(I)}+ \w^{(II)}|^2_{H^1}  \le | \w^{(I)} |^2_{H^1}+| \w^{(II)} |^2_{H^1} \leq
\frac{C}{\beta^2} \left< A_{\Delta\Delta}\w^{(I)}_\Delta,\w^{(I)}_\Delta\right> + \frac{C}{\beta^2 h^2}\left<g_{p_\Gamma},g_{p_\Gamma}\right>,
\]
where we used \EQ{uOne}, \EQ{uIbound}, \EQ{uuOne}, and \EQ{uuTwo} for the last inequality.

On the other hand, we have
\begin{eqnarray*}
\left< M^{-1}y,y \right> & = &
\frac{1}{h^2}\left<g_{p_\Gamma}, g_{p_\Gamma}\right>+ g_\lambda^T
M^{-1}_{1,\lambda}g_\lambda
= \frac{1}{h^2}\left<g_{p_\Gamma},g_{p_\Gamma}\right>+g_{\lambda}^T
B_{\Delta,D}\Rtilde_{\Delta}\Atilde\Rtilde_{\Delta}^TB_{\Delta,D}^Tg_{\lambda}\\
& = & \frac{1}{h^2}\left<g_{p_\Gamma},g_{p_\Gamma}\right> + \left<A_{\Delta\Delta} \vvec{w}^{(I)}_\Delta, \vvec{w}^{(I)}_\Delta\right>. \qquad \Box
\end{eqnarray*}

We also need the following lemma.
\begin{mylemma}
\label{lemma:BcW0}
For any $w = \left( {\bf w}_I, ~p_I, ~{\bf w}_{\Delta}, ~{\bf w}_{\Pi} \right) \in \Wtilde_0$, $B_C w \in R_G$.
\end{mylemma}

\beginproof
We know for any $\left( {\bf f}_I, ~{\bf f}_{\Delta}, ~{\bf f}_{\Pi} \right) \in \vvec{W}_I \bigoplus \vvec{W}_\Delta \bigoplus \vvec{W}_\Pi$, $g$ defined by \EQ{gvec} is in $R_G$.
For any $w = \left( {\bf w}_I, ~p_I, ~{\bf w}_{\Delta}, ~{\bf w}_{\Pi} \right) \in \Wtilde_0$, from the definition of $\Atilde$ in \EQ{AtildeBc}, there always exists $\left( {\bf f}_I, ~{\bf f}_{\Delta}, ~{\bf f}_{\Pi} \right) \in \vvec{W}_I \bigoplus \vvec{W}_\Delta \bigoplus \vvec{W}_\Pi$, such that
\[
\Atilde w = \left[
\begin{array}{l}
{\bf f}_I        \\[0.8ex]
0                \\[0.8ex]
{\bf f}_{\Delta} \\[0.8ex]
{\bf f}_\Pi
\end{array}
\right], \quad \mbox{i.e.,} \quad w = \Atilde^{-1} \left[
\begin{array}{l}
{\bf f}_I        \\[0.8ex]
0                \\[0.8ex]
{\bf f}_{\Delta} \\[0.8ex]
{\bf f}_\Pi
\end{array}
\right].
\]
Taking such $\left( {\bf f}_I, ~{\bf f}_{\Delta}, ~{\bf f}_{\Pi} \right)$, $g$ defined in \EQ{gvec} is $B_C w$.
$\qquad \Box$

The following lemma is an immediate result of Lemmas \LA{lower} and \LA{BcW0}.

\begin{mylemma}
\label{lemma:RGBcW0} The space $R_G$ is the same as the range of $B_C$ applied on $\Wtilde_0$.
\end{mylemma}

The condition number bound of the preconditioned operator $M^{-1} G$ is given in the following theorem.
\begin{mytheorem}
\label{theorem:tcond} For all $x = (p_{\Gamma}, ~\lambda) \in R_{M^{-1} G}$,
\[
C \beta^2 \left<Mx,x \right>\leq \left< G x,x \right> \leq \Phi(H/h) \left<
Mx,x \right>,
\]
where $\Phi(H/h)$ is as defined in Lemma \LA{jump}, $\beta$ as in \EQ{infsupMatrix}.
\end{mytheorem}

\beginproof
\[
\left< Gx,x\right> = x^T B_C\Atilde^{-1}B_C^Tx= x^T B_C \Atilde^{-1} \Atilde \Atilde^{-1} B_C^Tx =
\left< \Atilde^{-1} B_C^Tx, \Atilde^{-1} B_C^Tx\right>_{\Atilde}.
\]
Since $\Atilde^{-1} B_C^Tx \in \Wtilde_0$ and $\left< \cdot, \cdot \right>_{\Atilde}$ defines an inner product on $\Wtilde_0$, we have
\be \label{equation:Gnorm}
\left< Gx,x\right> = \max_{v \in \Wtilde_0, v \neq 0} \frac{\left<v, \Atilde^{-1} B_C^Tx \right>^2_{\widetilde{A}}}{\left<v, v \right>_{\widetilde{A}}} =\max_{v \in \Wtilde_0, v \neq 0}
\frac{\left<B_Cv,x\right>^2}{\left<\Atilde v,v\right>} .
\ee

{\it Lower bound:} From Lemma \ref{lemma:lower}, we know that for any given $y = (g_{p_{\Gamma}}, g_\lambda) \in R_G$, there exits
$w \in \Wtilde_0$, such that $B_C w = y$ and $\left <\Atilde w, w\right> \le \frac{C}{\beta^2}  \left< M^{-1}y, y \right>$. From  \EQ{Gnorm}, we have
$$
\left< Gx,x\right> \ge\frac{\left<B_Cw,x\right>^2}{\left<\Atilde w,w\right>}
\ge C \beta^2 \frac{\left<y,x\right>^2}{\left<M^{-1}y,y\right>}.
$$
Since $y$ is arbitrary, using Lemma \ref{lemma:m1RG}, we have
$$\left< Gx,x\right> \ge C  \beta^2 \max_{y \in R_G, y \neq 0} \frac{\left<y,x\right>^2}{\left<M^{-1}y,y\right>} = C  \beta^2 \left<Mx,x \right>.$$

{\it Upper bound:} From \EQ{Gnorm}, Lemmas \LA{upper}, \LA{RGBcW0}, and \LA{m1RG},  we have
\begin{eqnarray*}
\left< Gx,x\right> & \le & \Phi(H/h)\max_{v\in \Wtilde_0, v \neq 0} \frac{\left<B_Cv,x\right>^2}{\left<M^{-1} B_Cv,B_Cv\right>} \\
& = & \Phi(H/h) \max_{y\in R_G, y \neq 0} \frac{\left<y,x\right>^2}{\left<M^{-1} y,y\right>} = \Phi(H/h) \left<Mx,x\right>. \qquad \Box
\end{eqnarray*}

\begin{myremark}
\label{remark:rates}
From Theorem \ref{theorem:tcond} and Remark \ref{remark:Phi}, we can
see that the condition number bound of the preconditioned operator
$M^{-1} G$ is independent of the number of subdomains when $H/h$ is
fixed. If only subdomain corner velocities are chosen as
coarse level primal variables in the algorithm, the upper eigenvalue
bound of the preconditioned operator depends on $H/h$ in terms of
$(H/h) (1 + \log{(H/h)})$; if  both subdomain corner and edge-average velocity degrees of freedom are chosen as primal variables, the upper eigenvalue bound grows as $H/h$.
\end{myremark}

\begin{myremark}
\label{remark:discon}
With only minor modifications, the algorithm proposed in this paper and its analysis apply equally well to the discontinuous pressure case. In that situation, $p_\Gamma$ and the blocks related to it in \EQ{bigeq} can simply be replaced by the vector containing subdomain constant pressures and its corresponding blocks, respectively. The formulation of the algorithm then follows the same way as presented in Section \ref{section:Gmatrix}, and the same condition number bounds as in Theorem \ref{theorem:tcond} will be obtained. Numerical experiments of our algorithm for the discontinuous pressure case will also be reported in the next section.
\end{myremark}

\begin{myremark}
\label{remark:kim}
The same condition number bound has been proved by Kim and Lee~\cite[with Park]{kim102, kim10} for their FETI-DP algorithms for solving incompressible Stokes equations. In their algorithms, discontinuous pressure is considered and their approaches do not apply to the continuous pressure case.
\end{myremark}

\begin{myremark}
\label{remark:coarse}
We also note that, no additional coarse level degrees of freedom, except those necessary for solving positive definite elliptic problems, are required in our algorithm to achieve a scalable convergence rate. For example, for two-dimensional problems, it is sufficient to include only the subdomain corner velocity degrees of freedom in the coarse level problem. This represents a progress compared with earlier work, e.g., \cite{li05, li06}, where additional continuity constraints enforcing the  divergence-free conditions on subdomain boundaries are required in the coarse level problem. Reduction in the coarse level problem size has also been achieved for algorithms discussed in \cite{Doh09, Doh10, kim11, kim102, kim10}, even though discontinuous pressures are considered there.
\end{myremark}

\section{Numerical experiments}
\label{section:numerics}

We consider solving the incompressible Stokes problem \EQ{Stokes} in the square domain $\Omega=[0,1]\times
[0,1]$. Zero Dirichlet boundary condition is used. The right-hand side function $\vvec{f}$ is chosen such that the exact solution is
$$\u=\left[\begin{array}{c}
\sin^3(\pi x)\sin^2(\pi y)\cos(\pi y)\\[0.8ex]
-\sin^2(\pi x)\sin^3(\pi y)\cos(\pi x)
\end{array}\right] \quad \mbox{and}\quad
p=x^2-y^2.
$$

The modified Taylor-Hood mixed finite element, as shown in Figure~\ref{figure:TaylorHood}, is used for the finite element solution. The preconditioned system \EQ{Mspd} is solved by the CG iteration; the iteration is stopped
when the $L^2-$norm of the residual is reduced by a factor of $10^{-6}$.

Table \ref{table:M1} shows the minimum and maximum eigenvalues of the iteration matrix $M^{-1} G$, and the iteration counts. The coarse level variable space in this experiment is spanned by the subdomain corner velocities. We can see from Table \ref{table:M1} that the minimum eigenvalue is independent of the mesh size. The maximum eigenvalue is independent of the number of subdomains for fixed $H/h$; for fixed number of subdomains, it depends on $H/h$, presumably in the order of $(H/h) (1 + \log{(H/h)})$ as predicted in Remark \ref{remark:rates}.

\begin{table}[t]
\caption{\label{table:M1}  Solving \EQ{Mspd}, with only subdomain corner velocities in coarse space.}
\centering
\begin{tabular}{ccccc}
\quad $H/h$ (fixed) \quad & \quad \#sub \quad  & \quad $\lambda_{min}$
& \quad $\lambda_{max}$ \quad & \quad iteration  \quad \\[1.2ex]
\hline \\
8  & $4  \times  4$ & 0.35 &  8.92 & 21 \\[1.2ex]
  & $8  \times  8$ & 0.35 & 10.07 & 28 \\[1.2ex]
  & $16 \times 16$ & 0.35 & 10.23 & 29 \\[1.2ex]
  & $24 \times 24$ & 0.35 & 10.30 & 29 \\[1.2ex]
  & $32 \times 32$ & 0.35 & 10.33 & 29 \\
\hline \\[1.2ex]
\quad \#sub (fixed) \quad & \quad $H/h$ \quad  & \quad $\lambda_{min}$
& \quad $\lambda_{max}$ \quad & \quad iteration \quad \\[1.2ex]
\hline \\
$8  \times  8$ &  4 & 0.30 &  4.22 & 21 \\[1.2ex]
              &  8 & 0.35 & 10.07 & 28 \\[1.2ex]
              & 16 & 0.35 & 24.22 & 36 \\[1.2ex]
              & 24 & 0.35 & 40.12 & 43 \\[1.2ex]
              & 32 & 0.35 & 57.15 & 50 \\
\hline \\
& & & & \\
\end{tabular}
\end{table}

\begin{table}[t]
\caption{\label{table:M2}  Solving \EQ{Mspd}, with both subdomain corner and edge-average velocities in coarse space.}
\centering
\begin{tabular}{ccccc}
\quad $H/h$ (fixed) \quad & \quad \#sub \quad  & \quad $\lambda_{min}$
& \quad $\lambda_{max}$ \quad & \quad iteration  \quad \\[1.2ex]
\hline \\
8  & $4  \times  4$ & 0.36 &  4.29 & 17 \\[1.2ex]
  & $8  \times  8$ & 0.36 & 5.29 & 21 \\[1.2ex]
  & $16 \times 16$ & 0.36 & 5.56 & 21 \\[1.2ex]
  & $24 \times 24$ & 0.36 & 5.61 & 21 \\[1.2ex]
  & $32 \times 32$ & 0.36 & 5.64 & 21 \\
\hline \\[1.2ex]
\quad \#sub (fixed) \quad & \quad $H/h$ \quad  & \quad $\lambda_{min}$
& \quad $\lambda_{max}$ \quad & \quad iteration \quad \\[1.2ex]
\hline \\
$8  \times  8$ &  4 & 0.33 &  4.00 & 18 \\[1.2ex]
              &  8 & 0.36 & 5.29 & 21\\[1.2ex]
              & 16 & 0.36 & 11.63 & 26 \\[1.2ex]
              & 24 & 0.36 & 18.67 & 31 \\[1.2ex]
              & 32 & 0.36 & 26.12 & 36 \\
\hline
\end{tabular}
\end{table}

For the experiment reported in Table \ref{table:M2}, the coarse level variable space is spanned by both the subdomain corner velocities and the subdomain edge-average velocity components. Even though the edge-average velocity components are not necessary for the analysis, including them in the coarse level problem improves the convergence rate, for which the maximum eigenvalue in Table \ref{table:M2} grows in the order of $H/h$, as discussed in Remark \ref{remark:rates}.

\begin{figure}[hbt]
\begin{center}
\includegraphics[scale=0.6]{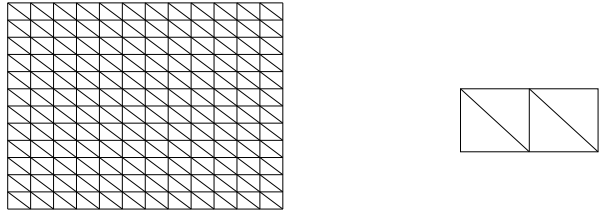}
\end{center}
\caption{The mesh and the mixed finite element.}
\label{fig:2dfem}
\end{figure}

Tables \ref{table:GDM1D} and \ref{table:GDM12D} show the performance of our algorithm for solving the same problem, but using a mixed finite element with discontinuous pressure. We use a uniform mesh of triangles, shown on the left in Figure~\ref{fig:2dfem}; the velocity finite element space contains the piecewise linear functions on the mesh and the pressure is a constant on each union of four triangles as shown on the right in the figure. The same mixed finite element has also been used in \cite{li06}.

Comparing Tables \ref{table:M1} and \ref{table:M2} with Tables \ref{table:GDM1D} and \ref{table:GDM12D}, we can see that the convergence rates of our algorithm, using either continuous or discontinuous pressure,  are quite similar.

\begin{table}[t]
\caption{\label{table:GDM1D}  Solving \EQ{Mspd}  (using discontinuous pressure), with only corner constraints.}
\centering
\begin{tabular}{ccccc}
\quad $H/h$ (fixed) \quad & \quad \#sub \quad  & \quad $\lambda_{min}$
& \quad $\lambda_{max}$ \quad & \quad iteration  \quad \\[1.2ex]
\hline \\
8  & $4  \times  4$ & 0.48 &  7.93 & 22 \\[1.2ex]
  & $8  \times  8$ & 0.48 & 9.00 & 25 \\[1.2ex]
  & $16 \times 16$ & 0.48 & 9.20 & 25 \\[1.2ex]
  & $24 \times 24$ & 0.48 & 9.20 & 25 \\[1.2ex]
  & $32 \times 32$ & 0.48 & 9.21 & 25 \\
\hline \\[1.2ex]
\quad \#sub (fixed) \quad & \quad $H/h$ \quad  & \quad $\lambda_{min}$
& \quad $\lambda_{max}$ \quad & \quad iteration \quad \\[1.2ex]
\hline \\
$8  \times  8$ &  4 & 0.41 &  3.91 & 19 \\[1.2ex]
              &  8 & 0.48 & 9.00 & 25 \\[1.2ex]
              & 16 & 0.49 & 21.39 & 36 \\[1.2ex]
              & 24 & 0.50 & 35.56 & 43 \\[1.2ex]
              & 32 & 0.50 & 50.87 & 50 \\
\hline
\end{tabular}
\end{table}

\begin{table}[h]
\caption{\label{table:GDM12D} Solving \EQ{Mspd}  (using discontinuous pressure),  with both corner and edge-average constraints.}
\centering
\begin{tabular}{ccccc}
\quad $H/h$ (fixed) \quad & \quad \#sub \quad  & \quad $\lambda_{min}$
& \quad $\lambda_{max}$ \quad & \quad iteration  \quad \\[1.2ex]
\hline \\
8  & $4  \times  4$ & 0.48 &  3.78 & 17 \\[1.2ex]
  & $8  \times  8$ & 0.49 & 4.47 & 18 \\[1.2ex]
  & $16 \times 16$ & 0.49 &4.68& 19 \\[1.2ex]
  & $24 \times 24$ & 0.50 & 4.77 & 19 \\[1.2ex]
  & $32 \times 32$ & 0.50 & 4.80 & 19\\
\hline \\[1.2ex]
\quad \#sub (fixed) \quad & \quad $H/h$ \quad  & \quad $\lambda_{min}$
& \quad $\lambda_{max}$ \quad & \quad iteration \quad \\[1.2ex]
\hline \\
$8  \times  8$ &  4 & 0.43 &  2.80 & 16 \\[1.2ex]
              &  8 & 0.49 & 4.47 & 18\\[1.2ex]
              & 16 & 0.50 & 9.85 & 26 \\[1.2ex]
              & 24 & 0.50 & 16.05 & 32 \\[1.2ex]
              & 32 & 0.50 & 22.67 & 37 \\
\hline
\end{tabular}
\end{table}

\section*{Acknowledgment} The authors are very grateful to Olof Widlund and Clark Dohrmann for their suggestion of this problem.

\end{document}